\def\lastrevised{April 2009.}
\documentclass[12pt]{article}
\usepackage{amssymb}

\def\header{The hierarchy of $\omega_1$-Borel sets}
\def\al{\alpha}
\def\be{\beta}
\def\bpi(#1){{\bf \Pi}^0_{#1}}
\def\bsi(#1){{\bf \Sigma}^0_{#1}}

\def\comp(#1){{\sim}{\kern -.01in}{#1}}
\def\concat#1{\hat{\phantom{a}}\la #1\ra}
\def\conj{\wedge\kern -.1in\wedge}
\def\continuum{{\mathfrak c}}
\def\de{\delta}

\def\elemsub{\preceq}
\def\fin(#1){{\rm Fn}(#1,2)}
\def\forces{{\;\Vdash}}
\def\ga{\gamma}
\def\ka{\kappa}
\def\la{\langle}
\def\om{\omega}
\def\poset{{\mathcal Q}}
\def\pow{{\mathcal P}}
\def\proof{\par\noindent Proof\par\noindent}
\def\pr{\prime}
\def\qed{\par\noindent QED\par}

\def\ra{\rangle}
\def\res{\upharpoonright}
\def\rmand{\mbox{ and }}
\def\rmiff{\mbox{ iff }}
\def\si{\sigma}
\def\sm{{\setminus}}
\def\starpi(#1){{\bf \Pi}^*_{#1}}
\def\starsi(#1){{\bf \Sigma}^*_{#1}}
\def\st{\;:\;} 
\def\su{\subseteq}

\newtheorem{theorem}{Theorem} 
\newtheorem{lemma}[theorem]{Lemma}
\newtheorem{question}[theorem]{Question}
\newtheorem{remark}[theorem]{Remark}

\begin{document}

\begin{center}
{\large The hierarchy of $\om_1$-Borel sets}
\end{center}

\begin{flushright}
Arnold W. Miller\footnote{
Thanks to the University of Florida Mathematics Department for
their support and especially Jindrich Zapletal, William Mitchell,
Jean A. Larson, and Douglas Cenzer for inviting me to the special
year in Logic 2006-07 during which most of this
work was done.
\par Last revised \lastrevised}
\end{flushright}

\def\address{\begin{flushleft}
Arnold W. Miller \\
miller@math.wisc.edu \\
http://www.math.wisc.edu/$\sim$miller\\
University of Wisconsin-Madison \\
Department of Mathematics, Van Vleck Hall \\
480 Lincoln Drive \\
Madison, Wisconsin 53706-1388 \\
\end{flushleft}}

\begin{center}
Abstract
\end{center}

\begin{quote}
We consider the $\om_1$-Borel subsets of the reals 
in models of ZFC.  This is the smallest
family of sets containing the open subsets of the $2^\om$ and
closed under $\om_1$ intersections  and $\om_1$ unions.  We show
that Martin's Axiom implies that the hierarchy of
$\om_1$-Borel sets has length $\om_2$.
We prove that in the Cohen real model the length of this
hierarchy is at least $\om_1$ but no more than $\om_1+1.$
\end{quote}

Some authors have considered $\om_1$-Borel sets in other spaces,
$\om_1^{\om_1}$ Mekler and Vaananen \cite{mekler}
and or completely metrizable spaces of uncountable density,
Willmott \cite{willmott}.  But in this paper we only consider
the space $2^\om$.

\bigskip\noindent
Define the levels of the
$\om_1$-Borel hierarchy of subsets of $2^\om$ as follows:

\begin{enumerate}
\item $\starsi(0)=\starpi(0)=$ clopen subsets of $2^\om$
\item $\starsi(\al)=\{\bigcup_{\be<\om_1}A_\be\st (A_\be:\be<\om_1)\in
(\starpi(<\al))^{\om_1}\}$
\item $\starpi(\al)=\{\bigcup_{\be<\om_1}A_\be\st (A_\be:\be<\om_1)\in
(\starsi(<\al))^{\om_1}\}$
\item  $\starpi(<\al)=\bigcup_{\be<\al}\starpi(\be)$
$\;\;\;\;\starsi(<\al)=\bigcup_{\be<\al}\starsi(\be)$
\end{enumerate}

\noindent The length of this hierarchy is
the smallest $\al\geq 1$ such that
$$\starpi(\al)=\starsi(\al).$$
It is easy to show that if
$\al<\om_2$ and every $\om_1$-Borel set is $\starpi(<\al)$, then
$\starpi(\be)=\starsi(\be)$ for some $\be<\al$, i.e., bounded
hierarchies must have a top class
(see Miller \cite{onthelength} Proposition 4 p.235).

The classes
$\starpi(1)$ and $\starsi(1)$ are the ordinary closed sets and open sets,
respectively, so the length of the hierarchy of $\om_1$-Borel sets is at
least $2$.

Assuming the continuum hypothesis,
$\starpi(2)=\starsi(2)=\pow(2^\om)$, so CH implies
the order of the hierarchy is $2$.
It also known
to be consistent that
$$\starpi(3)=\starsi(3)=\pow(2^\om)\rmand \starpi(2)\neq\starsi(2)$$
see Steprans \cite{step}.  In Stepran's model,
the continuum is $\aleph_{\om_1}$.   Carlson \cite{carlson} showed
that if subset of $2^\om$ is $\om_1$-Borel, then the cofinality
of the continuum must be $\om_1$.  Stepran's model was used
earlier by Bukovsky \cite{buk} and latter by
Miller-Prikry \cite{milprik}.

The following is an open question from
Brendle, Larson, and
Todorcevic \cite{larson}.

\begin{question}
Is it consistent with the negation of the continuum hypothesis
that $\starpi(2)=\starsi(2)$?
\end{question}

Steprans noted that it would be too much to ask for
$$\neg CH\;\;+\;\;\starpi(2)=\starsi(2)=\pow(2^\om)$$
since a $\starsi(2)$ set, i.e., an $\om_1$ union of
closed sets, of size greater than $\om_1$ would have to
contain a perfect subset, hence $\neg CH$ implies a
Bernstein set cannot be $\starsi(2)$.
It is also known that $\starpi(2)\neq\starsi(2)$ in
the iterated Sacks model, see Ciesielski and Pawlikowski \cite{cp}.

\begin{theorem}\label{ma-hier}
(MA$_{\om_1}$)  $\starpi(\al)\neq\starsi(\al)$ for every $\al<\om_2$.
\end{theorem}

We prove this using the following two lemmas.  A well-known consequence
of MA$_{\om_1}$ is that every subset $Q\su 2^\om$
of size $\om_1$ is a
Q-set, i.e., for every subset $X\su Q$ there is a
$G_\de$ set $G\su 2^\om$ with $G\cap Q=X$ (see Fleissner and Miller \cite{qsets}).

\begin{lemma}\label{qset} Suppose there exists a Q-set
of size $\om_1$. Then
there exists an onto map $F:2^\om\to 2^{\om_1}$ such
for every subbasic clopen set $C\su 2^{\om_1}$
the set $F^{-1}(C)$ is either $G_\de$ or $F_\si$.
\end{lemma}
\proof
Fix $Q=\{u_\al\in 2^\om\st\al<\om_1\}$ a Q-set.
Let $G\su 2^\om\times 2^\om$ be a universal
$G_\de$ set, i.e., $G$ is $G_\de$ and for every $G_\de$ set
$H\su 2^\om$ there exists $x\in 2^\om$ with $G_x=H$.
Define $F$ as follows, given $x\in 2^\om$ let
$$F(x)(\alpha)=1 \rmiff  u_\al\in G_x $$
If $C$ is a subbasic clopen set, then for some $\al$ and $i=0$ or $i=1$
$$C_{\al,i}=\{p\in 2^{\om_1}\st p(\al)=i\}.$$
Then for $i=1$
$$F^{-1}(C_{\al,1})=\{x\st u_\al\in G_x \}$$
which is a $G_\de$  set.  
Since $C_{\al,0}$ is the complement of $C_{\al,1}$ we have that
$F^{-1}(C_{\al,0})$ is an $F_{\si}$-set

Finally, we note that since $Q$ is a Q-set, i.e., every
subset is a relative $G_\de$, it follows that $F$ is onto.
\qed

The next Lemma is true without any additional
assumptions beyond ZFC.  Its proof is a generalization
of Lebesgue's 1905 proof
(see Kechris \cite{kechris} p.168) for the standard Borel hierarchy.

\begin{lemma}\label{univ}
For any $\al$ with $0<\al<\om_2$ there exists a $\starsi(\al)$ set
$U\su 2^{\om_1}\times 2^\om$ which is universal for
$\starsi(\al)$ subsets of $2^\om$, i.e., for any $Q\su 2^\om$ which
is $\starsi(\al)$ there exists $p\in 2^{\om_1}$ with
$U_p=Q$.  Similarly, there is a universal $\starpi(\al)$ set.
\end{lemma}
\proof
The proof is by induction on $\al$.  Note that the complement
of a universal $\starsi(\al)$ set is a universal $\starpi(\al)$-set.

For $\al=1$, $\starsi(\al)$ is just the open sets. There is a
universal open set $V\su 2^\om\times 2^\om$.  Put
$$U=\{(p,x)\in 2^{\om_1}\times 2^\om\st (p\res \om,x)\in V\}$$

For $\al$ such that $2\leq \al<\om_2$ proceed as follows.
Let $(\de_\be<\al:\be<\om_1)$ have
the property that for every $\ga<\al$ there are $\om_1$  many
$\de_\be \geq \ga$.  It follows that for every $\starsi(\al)$
set $Q\su 2^\om$ there is
$(Q_\be\in \starpi(\de_\be)\st \be<\om_1)$ with
$$Q=\bigcup_{\be<\om_1}Q_\be.$$
By induction, there are
$U_\be\su 2^{\om_1}\times 2^\om$ universal $\starpi(\de_\be)$
sets.  Let $a:\om_1\times\om_1\to\om_1$ be a bijection.
For each $\be$ define
$$\pi_\be: 2^{\om_1}\times 2^\om\to  2^{\om_1}\times 2^\om,\;\;
(p,x)\mapsto (q,x)$$
where $q(\al)=p(a(\be,\al))$.
Put
$$U=\bigcup_{\be<\om_1}\pi_\be^{-1}(U_\be)$$
then $U$ will be a universal $\starsi(\al)$ set.
\qed

Now we prove Theorem \ref{ma-hier}.  Suppose for contradiction, that
every $\om_1$-Borel set is $\starsi(\al)$ for some fixed
$\al<\om_2$.  Let $U\su 2^{\om_1}\times 2^\om$ be a universal
$\starsi(\al)$ and define
$$V=\{(x,y)\in 2^\om\times 2^\om\st (F(x),y)\in U\}.$$
Then $V$ is an $\om_1$-Borel set (although not necessarily at
the $\starsi(\al)$) because the preimage of any
clopen box $C\times D$ is $\om_1$-Borel by Lemma \ref{qset}.
Define
$$D=\{x:(x,x)\notin V\}.$$
But then $D$ is $\om_1$-Borel but not $\starsi(\al)$.  We see
this by the usual
diagonal argument that if $D=U_p$, then since $F$ is onto
there would be $x\in 2^\om$ such that $F(x)=p$ but then
$$x\in D \rmiff (F(x),x)\notin U \rmiff x\notin U_p
\rmiff x\notin D.$$
\qed

\begin{remark}
Note that in the proof $V\su 2^\om\times 2^\om$ is a
$\starsi(2+\al)$-set, since the preimage of a clopen set
under $F$ is ${\bf\Delta}_3^0$.  Hence for levels $\al\geq\om$
the set $V$ is a $\starsi(\al)$ set which is universal 
for $\starsi(\al)$ sets.
\end{remark}

\begin{remark}
Our result easily generalizes to show that MA implies
that for any $\ka$ a cardinal with $\om\leq \ka <|2^\om|$ the
$\ka$-Borel hierarchy has length $\ka^+$.  This implies
that for any $\ka_1<\ka_2$ there are $\ka_2$-Borel sets
which are not $\ka_1$-Borel.\footnote{Since $\ka_2$-Borel sets at level
$\ka_1^+$ or higher cannot be $\ka_1$-Borel.}
It is also true for the Cohen real model that
for $\om\leq \ka_1<\ka_2<|2^\om|$ that there are $\ka_2$-Borel sets
which are not $\ka_1$-Borel.
\end{remark}

\begin{question}
Suppose MA and the continuum, $\continuum=|2^\om|$,
is a weakly inaccessible cardinal.  What is the
length\footnote{The argument of Lemma \ref{lowerbound} shows
that it is at least $\om_1$.} of
the hierarchy of ($<\continuum$)-Borel sets?
\end{question}

\begin{theorem} \label{cohenmodel}
In the Cohen real model every $\om_1$-Borel set is
$\starsi(\om_1+1)$ and there is a $\starsi(\om_1)$ set which
is not in $\starsi(<\om_1)$.
\end{theorem}
\proof

We state the lower bound separately as Theorem \ref{lowerbound}.

We will use Steel forcing with tagged trees (Steel \cite{steel})
similarly to its use in Stern \cite{stern}.  Stern proved that
assuming MA$_{\om_1}$ an $\om_1$ union of
$\bsi(\al)$ sets which is Borel, must be $\bsi(\al)$.
Since Steel forcing is countable, he only really needed MA(ctble).
Similar results are proved in Solecki \cite{solecki} Cor 2.3 and
Becker and Dougherty \cite{becker} Thm 2.
These authors do not consider
$\om_1$-Borel sets but are interested only in $\om_1$-unions
of ordinary Borel sets.

$MA_{\om_1}(ctble)$ stands for Martin's axiom for countable posets.
It says that for any countable poset and $\om_1$-family of dense sets
there is a filter meeting all the dense sets in the family.  
It is equivalent to saying that the real line cannot 
be covered by $\om_1$ nowhere
dense sets, see for example, Bartoszynski and Judah \cite{barto} p. 138.
It holds in any generic extension obtained with a finite 
support ccc iteration of cofinality at least $\om_2$.

\begin{theorem}\label{lowerbound}
Suppose $MA_{\om_1}(ctble)$ holds. Then for
any $\al<\om_1$ there is an ordinary Borel set
which is not $\starsi(\al)$.
\end{theorem}
\proof
We use Steel forcing with tagged trees\footnote{Sami \cite{sami} gives
a proof of Harrington's Theorem which does not use Steel forcing.}
similarly to the way it is described in
Harrington \cite{harrington}.

For any countable ordinal $\al$ define
$\poset(\al)$ to be the following countable poset.
Elements of $\poset(\al)$ have the form $(t,h)$ where
$t$ is a finite subtree of $\om^{<\om}$ and $h:t\to\al\cup\{\infty\}$
is called a tagging.  The ordering on $\al\cup\{\infty\}$ is
$\infty<\infty$ and $\be<\infty$ for each ordinal $\be$ along
with the usual ordering on pairs of ordinals from $\al$.  A tagging
$h$ is a rank function which means it satisfies:
if $\si,\tau\in t$ and $\si$ is a strict initial segment of $\tau$,
then $h(\si)>h(\tau)$.\footnote{We differ from \cite{harrington}
by not requiring that $h(\la\ra)=\infty$.}

The ordering on $\poset(\alpha)$ is
$p\leq q$ ($p$ extends $q$) iff
\begin{enumerate}
\item $t_q\su t_p$ and
\item $h_q\su h_p$.
\end{enumerate}
Note that nodes tagged with $\infty$ can always be extended and
tagged with $\infty$ or any element of $\al$.
\footnote{Harrington \cite{harrington} makes the additional requirement
that the top node, $\la\ra$, be tagged with $\infty$, but this
is unnecessary and makes our proof clumsy, as in
Miller \cite{proj} Lemma 4.4.}

Now suppose that $G$ is $\poset(\al)$ generic
over $M$.  Define
\begin{enumerate}
\item $T=T_G=\{\si\st \exists (t,h)\in G \;\; \si\in t\}$
\item $H=H_G:T_G\to\al\cup\{\infty\}$ by $H(\si)=h(\si)$ for
any $h$ such that there exists $(t,h)\in G$ with $\si\in t$.
\end{enumerate}

It is easily seen by a density argument that $H$ is a rank function on the
tree $T$ where the symbol $\infty$ gets attached to the nodes of $T$
which can be extended to an infinite branch.

Define $p(\beta)$ for $\beta\leq\alpha$ and $p\in \poset(\alpha)$
by $p(\beta)=(t,h_\beta)$ where $p=(t,h)$ and
$$h_\beta(s)=\left\{
\begin{array}{ll}
h(s) & \mbox{if } h(s)<\om \cdot\beta \\
\infty & \mbox{otherwise}
\end{array}
\right.$$

\begin{lemma}
(Retagging Lemma) Suppose $p_1,p_2\in \poset(\al)$ and
$\be+1\leq \al$ and $p_1(\be+1)=p_2(\be+1)$.  Then for
every $q_1\leq p_1$ there is $q_2\leq p_2$ such that
$q_1(\be)=q_2(\be)$.
\end{lemma}
\proof
Let $p_i=(t,h_i)$ for $i=1,2$ and suppose $q_1=(t^\pr,f_1)$.
We define $q_2=(t^\pr,f_2)$ as follows.
Put $f_2\res t=h_2$.  Fix $N<\om$ greater than the height
of $t^\pr$.  For each $\si\in t^\pr\sm t$ let $\tau\su \si$ be
the longest initial segment of $\si$ which is in $t$.

Case 1.
If $h_1(\tau)<\om(\be+1)$, then by assumption, $h_2(\tau)=h_1(\tau)$
and we can define $f_2(\si)=f_1(\si)$.

Case 2.
$h_1(\tau)\geq \om(\be+1)$, then by assumption,  $h_2(\tau)\geq \om(\be+1)$.

(a) If $h_1(\si)<\om\be$, then we put $h_2(\si)=h_1(\si)$.

(b) Otherwise $\om\be\leq h_1(\si)$ and
we put $h_2(\si)=\om\be+(N-|\si|)$.  Note that in this
case when we look at $q_i(\be)$ these $\si$ will be retagged
with $\infty$.
\qed

Fix $\al<\om_1$ and let $ T $ be the usual $\poset(\al)$-name
for the generic
tree $T_G$:
$$ T =\{(p,\check{s}):s\in t_p\mbox{ where } p=(t_p,h_p)\in
\poset(\al)\}.$$

The following is the main property of Steel forcing.
We identify $\pow(\om^{<\om})$ with $2^\om$.

\begin{lemma}\label{mainsteel}
Suppose $p,q\in\poset(\al)$,
$1+\be\leq\al$, $p(1+\be)=q(1+\be)$, and
$B\su \pow(\om^{<\om})$ is $\starpi(\be)$ set coded
in the ground model.\footnote{There are many ways to
code Borel (or more generally $\ka$-Borel) sets.
Solovay \cite{solovay} p.25 gives
a clear definition of coding and absoluteness which is similar
to what we use in the proof of Lemma \ref{univ}. Harrington \cite{harrington}
Definition 2.5 and Steel \cite{steel} code using infinitary
propositional logic. We like to use well-founded trees
as in Lemma \ref{upperbound}.}
Then
$$p\forces T \in B \rmiff
q\forces T \in B .$$
\end{lemma}
\proof
This is proved by induction on $\be$.

For $\be=0$ we
take for $\starpi(0)$ basic clopen subsets of
$\pow(\om^{<\om})$.  This means that for some
pair $F_0,F_1$ of disjoint finite subsets of $\om^{<\om}$ that
$$B=\{X\su \om^{<\om} \;:\; F_0\su X \rmand F_1\cap X=\emptyset\}.$$
So the statement $X\in B$ is a finite conjunction of statements
of the form $\si\in X$ or $\si\notin X$.  But note that:

\begin{enumerate}
\item $p\forces \check{\si}\in T $
iff $\si=\la\ra$, $\si\in t_p$, or
$\tau\in t_p$ and $h_p(\tau)>0$ where $\tau$ is the initial segment
of $\si$ of length exactly one less than $\si$.
\item $p\forces \check{\si}\notin T $
iff there exists $\tau\su\si$ with $\tau\in t_p$ and
$h_p(\tau)<|\si|-|\tau|$.
\end{enumerate}

Both of these are preserved when we look at $p(1)$.
Hence if $p(1)=q(1)$ then
$$p\forces  T \in B \rmiff
q\forces  T \in B .$$

For $\be>0$ suppose that $B$ is $\starpi(\be+1)$
and coded in the ground model. Working in the ground model
let $B=\bigcap_{\al<\om_1}\comp(B)_\al$
where\footnote{We use $\comp(B)$ to denote the complement of $B$.}
each $B_\al$ is $\starpi(<1+\be)$.  And suppose for contradiction
that $p_2\forces  T \in B$
but $p_1$ does not force this.  Then there exists a $q_1\leq p_1$
and $\al<\om_1$ such that
$$q_1\forces  T \in B_\al.$$
And suppose that $B_\al$ is $\starpi(1+\ga)$ where
$\ga<\be$. Since $1+\ga+1\leq 1+\be$,
by the retagging lemma we may find $q_2\leq p_2$ with
$q_1(1+\ga)=q_2(1+\ga)$.  By inductive hypothesis
$$q_2\forces  T \in B_\al$$
which contradicts that
$$p_2\forces  T \in B\supseteq \comp(B)_\al.$$
\qed

Suppose for contradiction that in
the Cohen real model there is an $\al_0<\om_1$ such that
every $\om_1$-Borel set
is $\starpi(\al_0)$.  It well-known that for every countable ordinal $\al$
the set
$$WF_{\al}=\{T\su\om^{<\om}: T \mbox{ is a well-founded tree of rank }
\al\}$$
is an (ordinary) Borel set.\footnote{The exact Borel
class is computed in Stern \cite{stern} and Miller \cite{orbits}.}
Consequently it must be a
$\starpi(\al_0)$-set.
Fix a countable $\al>\al_0\cdot\om$.  Take a sufficiently large\footnote{
For example $\ka=\beth_\om^+$.}
regular cardinal
$\ka$ and let $H_\ka$ be the sets whose transitive closure
has cardinality less than $\ka$.
Take $N$ to be
an elementary
substructure of $V_\ka$ of cardinality $\om_1$ which
contains $\al+1$. Then $N$ will contain a code for $B$ the $\starpi(\al_0)$
set $WF_{\al}$.
Let $M$ be the transitive collapse of $N$ and
consider forcing over $M$ with $\poset(\al+1)$.  Since we are
assuming MA(ctbl), for any $p\in\poset(\al+1)$ there
is a $G$ $\poset(\al+1)$-generic over the ground model $M$
with $p\in G$.  So take such a $G$ with $H_G(\la\ra)=\al$.
Then $T_G$ is a well-founded tree of rank $\al$ and
so $T_G\in WF_{\al}$.  By absoluteness
$$M[G]\models T_G\in B$$
and so there must be a $p\in G$ such that
$$p\forces  T \in B.$$
But consider $q=p(\al)$.  Note that $h_q(\la\ra)=\infty$.
Consequently, for any $G^\pr$ which is $\poset(\al+1)$-generic over $M$
with $q\in G^\pr$, the tree $T_{G^\pr}$ is not even well-founded
and hence
$$M[G^\pr]\models T_{G^\pr}\notin B.$$
But this means that
$$q\forces  T \notin B$$
which contradicts Lemma \ref{mainsteel}.
\qed

Next we prove an upper bound on the $\om_1$-Borel hierarchy
in the Cohen real model.  Our argument uses some ideas
employed by Carlson \cite{carlson}.

\begin{lemma} \label{upperbound}
In the Cohen real model for any $\om_1$-Borel set $B$ there exists
$\om_1$  ordinary Borel sets, $(B_\be:\be<\om_1)$,
such that $B$ is their limit:
$$B=\bigcup_{\al<\om_1}\bigcap_{\be>\al}B_\be=
\bigcap_{\al<\om_1}\bigcup_{\be>\al}B_\be
$$
\end{lemma}
\proof
Let $B$ be coded by a well-founded tree $T\su \om_1^{<\om_1}$
with basic clopen sets $(s_\si\in 2^{<\om}:\si\in T^*)$ where
$T^*$ are the terminal nodes (or leaf nodes) of the tree $T$.
Then $T,(s_\si:\si\in T^*)$ codes $B$ as follows.
Define
$$B(\si)=[s_\si]=\{x\in 2^\om \st s_\si\su x\}$$
for $\si\in T^*$.  Then for nonterminal nodes of $T$ define
$$B(\si)=
\bigcap\{\comp(B(\si\concat\al)):\al<\om_1\rmand \si\concat\al\in T\}.$$
Finally, put $B=B({\la\ra})$.

Fix such a $T$ for $B$ and for any $\al<\om_1$ define
$B_\al$ inductively just as above but for the countable tree
$T\cap\al^{<\om}$.

We will show that for some closed unbounded set $C\su\om_1$
that $B$ is the $\om_1$-limit of $(B_\be:\be\in C)$.

By the Cohen real model we mean an model obtained by forcing
with $\fin(\om_2)$, the finite partial maps from
$\om_2$ into $2$, over a model of ZFC+GCH.  By standard arguments
using the countable chain condition and product Lemma, we
may without loss of generality assume that our code
for $B$, $T,(s_\si:\si\in T^*)$, is in the ground model $M$
a model of ZFC+GCH.
For any $x\in M[G]\cap 2^\om$ (where $G$ is $\fin(\om_2)$-generic over $M$
there is an $H\in M[G]$ which is $\fin(\om)$-generic over $M$ and
$x\in M[H]$.

Since the ground model $M$ satisfies CH,
there is a set of canonical names, CN, for elements of
$2^\om$ in the extension $M[H]$ has size $\om_1$.

Working in the ground model $M$ construct an continuous chain
$(N_\al:\al<\om_1)$
of countable elementary submodels of $H_{\om_2}$, with
the code for $B$, $T,(s_\si:\si\in T^*)$, in $N_0$,
$N_\al\elemsub N_\be$ and $N_\al\in N_\be$ for $\al<\be<\om_1$.
Note that it is automatically the case that every canonical name
is in some $N_\al$.

Now take for our club $C$ the set
$$C=\{\om_1\cap N_\al\st \al<\om_1\}.$$

Suppose that $x=\tau^H$ where $\tau\in N_\al$ and
$H$ is $\fin(\om)$-generic over $M$.  Let $M_\al$ be
the transitive collapse of $N_\al$.  By standard arguments
$H$ is $\fin(\om)$-generic over $M_\al$.
Note that ordinal $\delta=N_\cap\om_1$ is the
$\om_1$ of $M_\al$ i.e.,
$$M_\al\models \delta=\om_1.$$
Let
$p\in\fin(\om)$ be such that either
$$M_\al\models p\forces \tau\in B$$
or
$$M_\al\models p\forces \tau\in \comp(B).$$
Assume the former.  Note that
$B^{M_\al[H]}$=$B_\de\cap M[H]$.  And since it is forced
it must be that $x=\tau^H\in B_\de$.

For every $\be>\al$ the model $N_\be$ elementary superstructure of
$N_\al$ and hence that
$$M_\be\models p\forces \tau\in B$$
and for the same reason $x\in B_{\de^\pr}$ where
$\de^\pr$ is the $\om_1$ of $M_\be$.
\qed

\begin{remark}
Lemma \ref{upperbound} easily generalizes to the $\om_1$-Borel hierarchy
giving that every $\om_2$-Borel set is the $\om_2$ limit of $\om_1$-Borel
sets, and since each of them is at level $\om_1+1$, we get
an upperbound of $\om_1+2$ for the length of the $\om_2$-Borel hierarchy.
\end{remark}

\begin{remark}
Lemma \ref{upperbound} is also true in the random real model.
\end{remark}

\begin{remark}
In Steprans \cite{step} the hierarchy on the $\om_1$-Borel
sets is defined by letting the bottom level,
$\Pi^{\aleph_1}_0=\Sigma^{\aleph_1}_0$, be the family of all
ordinary Borel sets.  Lemma \ref{upperbound} shows that every
$\om_1$-Borel set in the Cohen real model is
$\Pi^{\aleph_1}_2$ and hence $\Sigma^{\aleph_1}_2$.  It is
easy to see that in this model there are
$\Sigma^{\aleph_1}_1$  sets  which are not
$\Pi^{\aleph_1}_1$, for example, any nonmeager subset of $2^\om$ of
size $\om_1$.
\end{remark}

\begin{remark}
In Miller \cite{onthelength} Theorem 34 and 54, it is shown consistent
for any countable ordinal $\al_0\geq 2$ to have separable metric space
$X$ such that every subset of $X$ is Borel and the Borel hierarchy on $X$
has length exactly $\al_0$. It is easy to show that if the set
$X\su 2^\om$ has cardinality at least $\om_2$ that for
each $\be<\al_0$ the generic
${\bf\Pi}^0_\be$ sets produced are not $\starsi(\be)$ relative to $X$.
Hence these spaces have order $\al_0$ in the relativized $\om_1$-Borel
hierarchy.  If we replace the use of almost disjoint forcing in
Steprans model \cite{step} Definition 2,
by $\Pi_{\al_0}^0$-forcing from Miller \cite{onthelength} p. 236,
then we get a model of ZFC in which
every subset of $2^\om$ is $\om_1$-Borel and
the $\om_1$-Borel hierarchy has length
at least $\al_0$ but no more than $\al_0+1$.  Similarly if
we change the Steprans model
by using $\Pi_\al^0$-forcing in the $\al$
model,
then in the resulting model every subset of $2^\om$ is $\om_1$-Borel and
the $\om_1$-Borel hierarchy has length at least $\om_1$ but
no more than $\om_1+1$.
\end{remark}

\begin{question}
Is possible to have a model of ZFC in which the $\om_1$-Borel hierarchy
has length $\al$ where $\om_1+2\leq \al<\om_2$?
\end{question}

\address

\end{document}